\documentclass[12pt]{amsart}
\usepackage{amsmath}
\usepackage{amsfonts,amssymb,amsthm, txfonts, pxfonts}
\def\struckint{\mathop{%
\def\mathpalette##1##2{\mathchoice{##1\displaystyle##2}%
  {##1\textstyle##2}{##1\scriptstyle##2}{##1\scriptscriptstyle##2}}%
\mathpalette
{\vbox\bgroup\baselineskip0pt\lineskiplimit-1000pt\lineskip-1000pt
\halign\bgroup\hfill$}
{##$\hfill\cr{\intop}\cr\diagup\cr\egroup\egroup}%
}\limits}

\newtheorem{theorem}{Theorem}
\newtheorem{lemma}[theorem]{Lemma}
\newtheorem{corollary}[theorem]{Corollary}
\newtheorem{conjecture}[theorem]{Conjecture}

\newtheorem{theorem-definition}[theorem]{Theorem-Definition}

\theoremstyle{remark}

\newtheorem{remark}[theorem]{Remark}

\DeclareMathOperator{\Sp}{Sp}
\DeclareMathOperator{\SL}{SL}
\DeclareMathOperator{\GL}{GL}
\DeclareMathOperator{\Gal}{Gal}

\begin{document}

%-------------- Author entries --------------------

\title[Polynomials, matrices, and automorphisms]{Counting reducible matrices, polynomials, and surface and free group automorphisms}

%Article title
%\shorttitle{Polynomials, matrices, and maps} % Shortened version for
                                             % headline title

\author{Igor Rivin}

\address{Department of Mathematics, Temple University, Philadelphia}

\email{rivin@math.temple.edu}

\thanks{I  would like to thank Ilya Kapovich for asking the
  questions on the irreducibility of  automorphisms of surfaces and
  free groups, and for suggesting that these questions may be
  fruitfully attacked by studying the action on homology.
  I would also like to thank Nick Katz for making him aware of Nick
Chavdarov's work, Sinai Robins for making him aware of Morris Newman's
classic book, and to Akshay Venkatesh and Peter Sarnak for enlightening
conversations. The author would also like to thank Benson Farb, Ilya Kapovich, and Lee Mosher for
comments on previous versions of this note.}

\curraddr{Mathematics Department, Princeton University, Princeton, NJ}

\date{\today}


%\subjclass{57M50, 32G15}

\begin{abstract}
We give upper bounds on the numbers of various classes of polynomials 
reducible over $\mathbb{Z}$ and over $\mathbb{Z}/{p \mathbb{Z}},$ and on the number of matrices in 
$\SL(n), \GL(n)$ and $\Sp(2n)$ with reducible characteristic
polynomials, and on polynomials with non-generic Galois groups. We use our result to show that a random (in the
appropriate sense) element of the mapping class group of a closed
surface is pseudo-Anosov, and that a random automorphism of a free
group is irreducible with irreducible powers\footnote{Following Handel and Mosher, we will call such an automorphism \emph{strongly irreducible}, thus avoiding the horrible acronym ``iwip''}
We also give a necessary condition
for all powers of an algebraic integers to be of the same degree, and
give a simple proof (in the Appendix) that the distribution of cycle
structures mod $p$ for polynomials with a restricted coefficient is
the same as that for general polynomials.
\end{abstract}

\maketitle

\section*{Introduction} 

In this paper we use simple algebraic, geometric, and probabilistic
ideas to investigate the probability that a random (in a suitable
sense) polynomial with integer coefficient is reducible (over
$\mathbb{Z}$) and that a random (in a suitable sense) matrix in one of
the classical groups $\GL(n, \mathbb{Z}), \SL(n, \mathbb{Z})$ or 
$\Sp(n, \mathbb{Z})$ and also in $M^{n\times n}(\mathbb{Z})$
has irreducible characteristic polynomial\footnote{To avoid encumbering the notation, we state the results for $\SL(n, \mathbb{Z}).$ The results and the proofs for $\GL(n, \mathbb{Z})$ are essentially identical}. We use these 
results (following an idea of I. Kapovich) to show that for 
generating set of the mapping class group, a sufficiently long random
product of generators is almost certainly pseudo-Anosov\footnote{A closely related result on the mapping class group was shown by completely different methods by J.~Maher in \cite{maheranosov}.}.

 The plan of the rest of the paper is as follows:
 In Section \ref{closure} we discuss some generalities on elimination 
 theory. In Section \ref{poly} we apply the results of Section 
 \ref{closure} to gain insight into sets of polynomials with factors of 
 certain types, in particular on the growth of the cardinality of 
 these sets as a function of height. In Section \ref{matrices} we apply the results of 
 Sections \ref{closure} and \ref{poly} to the study of sets of matrices 
 in $M^{n\times n}(\mathbb{Z})$ whose characteristic 
 polynomials are \emph{reducible}, and again, to get estimates on the 
 growth of these sets as a function of height (the size of 
 coeffcients). In Section \ref{matrices2}
we will use a quite different method to show that the density of 
 ``reducible'' elements in $\Sp(2n, \mathbb{Z})$ goes 
 to zero as a function of the \emph{combinatorial} distance of the elements to 
the identity. In Section \ref{mcg} we use our results to show that a
random element of the mapping class group of a closed surface of genus
$g$ is pseudo-Anosov. In Sectin \ref{strong} we show (using results of
the Appendix)
hat the characteristic polynomial of a
random matrix in $\GL(N, \mathbb{Z})$ has characteristic polynomial
with Galois group $S_N,$ and also that all powers of such a random
matrix have the same property. In Section \ref{free} we apply our
results to show that a random free group automorphism is strongly
irreducible (what is more commonly known in the trade as ``irreducible
with irreducible powers''.)

 \section{Generalities on elimination}
\label{closure}
 
Consider the following setup: we have a parametrized surface $S$ in $k^n$
(for $k$ an algebraically closed field), that is:
\begin{gather}
x_1 = f_1(s_1, \dotsc, s_m),\\
x_2 = f_2(s_1, \dotsc, s_m),\\
\vdots\\
x_n = f_n(s_1, \dotsc, s_m),
\end{gather}
where $f_1, \dotsc, f_n$ are polynomials in $s_1, \dotsc, s_m.$
It is reasonable to believe that $S$ is an algebraic $m$-dimensional
variety in $k^n,$ that is, the simultaneous zero-set of $n-m$
polynomial equations. That turns out to not be exactly true, but what
is true is that the \emph{Zariski closure} of $S$ is an
(at most) $m$-dimensional variety. For a proof of this \emph{Closure Theorem}
and plenty of examples see \cite{CLO_IVA}[Chapter 3].
 
 \section{Applications to polynomials}
\label{poly}
 Let $\mathcal{P}$ be the set of all monic  polynomials in one
variable of 
 degree $d$ over a field 
 $F,$ which have a polynomial factor with constant term $\alpha.$ 
Let us identify the set of all monic polynomials of degree $d$ with 
the affine space $F^{d}.$ Then, we have the following:
\begin{theorem}
    \label{hyper1}
    The set $\mathcal{P}$ is contained in an affine hypersurface of 
    $F^{d}.$
\end{theorem}

\begin{proof}
    Let \[p(x) = x^{d} + \sum_{i=0}^{d-1}a_{i} x^{i}¥ \in 
    \mathcal{P}.\] By assumption, $p(x) = q(x) r(x).$ Assume that the 
    degree of $q(x) = m,$ while the constant term of $q(x)$ equals 
    $\alpha.$ Writing \[q(x) = x^{m} + \sum_{j=1}^{m-1} b_{j}x^{j},\] 
    and \[r(x) = x^{d-m} + a_{0}/\alpha + \sum_{k=1}^{d-m-1}
c_{k}x^{k},\] we find ourselves exactly in the setting of Section
\ref{closure}. The proof is almost complete, except for the fact that
we do not know the degree of $q(x)$ \emph{a priori.} However, each
choice of $m$ gives us a polynomial $H_m$  vanishing at all the coefficient
sequences of reducible polynomials with a factor of degree $m,$   and
so the product of $H_m$ over all $m$ vanishes at \emph{all} the
coefficient sequences of reducible polynomials.
    \end{proof}

\section{Counting points on varieties}
\label{counting}

Let $S$ be a variety of dimension $m$ in $k^n.$.  Consider a reduction of $S$
modulo $p.$

\begin{theorem}[Lang-Weil, \cite{langweil}]
The number of $F_p$ points on $S$ grows as $O(p^m).$ The implied
constant is uniform (that is, it is a function of the dimension and
codimension of the variety only.
\end{theorem}

It should be noted that this gives an upper bound only. There might
well be \emph{no} $F_p$ points on $S.$

The following corollary is also classical (and easy):

\begin{corollary}
\label{heights}
Let $S$ be as above. Then the number of points of $S \cup
\mathbb{Z}^n$ all of whose coordinates do not exceed $B$ in absolute
value grows at most as $O(B^m).$
\end{corollary}

\begin{proof}
Pick $B.$ By Bertrand's postulate there is a prime $p,$ such that $4B
> p > 2 B.$ We know that every integer point of $S$ will give a
(distinct) point on the reduction of $S$ modulo $p$ (the converse, of
course, is not true). The result follows.
\end{proof}

We have used
\begin{theorem}[Bertrand's Postulate - proved by Chebyshev]
For any $N>3$ there exists at least one prime $p$ between $n$ and $2n-2.$
\end{theorem}

\section{More applications to polynomials}
\label{poly2}
The results in Section \ref{counting} combined with the results in
Section \ref{closure} immediately give the following results:

\begin{theorem}
\label{constone}
Let $P_1(d, B)$ be the set of polynomials of degree $d$ with integer coefficients
bounded in absolute value by $N$ and constant coefficient $1,$ and let
$R_1(d, B)$ be the set of polynomials reducible over $\mathbb{Z}$ with
the same coefficient bound. Then, $R_1$ lies on an algebraic
hypersurface $\mathbb{C}^{d-1}$ (where the coordinates are the
coefficients), and consequently
\[
\dfrac{R_1(B)}{P_1(B)} = O\left(\dfrac{1}{B}\right).
\]
\end{theorem}
\begin{proof}
A factor of a polynomial in $R_1(d, B)$ must have constant term
$\pm1,$ The statement now follows immediately from the results in
Sections \ref{closure} and \ref{poly}.
\end{proof}

\subsection{Arbitrary polynomials}
\label{arbitrarypoly}
What happens if we don't require the constant coefficient to be $1$?
Consider the set $F(d, a)$ of all monic polynomials of degree $d$ and
with constant term $a.$ Clearly, the constant term of a divisor of
such a polynomial must have constant term $d$ dividing $a,$ and so for
each $c | a$ we have a subvariety of $F(d, a)$ of polynomials having a
factor with constant term $c.$ The arguments above apply without
change, and we see that the number of such polynomials modulo $p$
grows at most as $O(p^{d-2}),$ where the constant is uniform. Denoting
the number of divisors of $a$ by $\tau(a),$ it is not hard to see that
$\tau(a) = o(a).$ Indeed, since the number of divisors is a
multiplicative function, 
\[\tau(n=p_1^{\alpha_1} \dots p_k^{\alpha_k}) = (\alpha_1 + 1) \dots
(\alpha_k + 1) < 2 \log_2 n,\] whereupin the assertion follows easily. 

So, it follows that for any $a,$ the set of reducible polynomials is a
union of $o(a)$ subvarieties of $F(d, a).$ To show that most
polynomials with coefficients bounded by $B$ in absolute value are
irreducible, we use Bertand's postulate to find a prime $p,$ such that
$2B<p <4B.$ Since the set of reducible polynomials lies on the union
of $o(B^2)$ codimension two subvarieties, their total number is
$o(p^d),$ while the total number of monic polynomials is $B^d \geq
2^{-d} p^d,$ and we have our result.

\section{Reciprocal Polynomials}
\label{recip}
We say that a polynomial $p(x) \in P_1(d)$ is  \emph{reciprocal} if 
$x^d p(1/x) = p(x)$ -- in other words, the list of coefficients of $p$
is the same read from left to right as from right to left. Reciprocal
polynomials can also be defined as follows:
A (monic) polynomial (of even degree $2n$) is reciprocal if it can be
written as
\[
\prod_{j=1}^n (x-r_i)(x-1/r_i) = \prod{j=1}^n (x^2 - (r_i + r_i^{-1})
x + 1).
\]

Notice that this means that \emph{every} recriprocal polynomial lies
on our ``factorization variety''\footnote{the author thanks N. Katz
for the suggestion of using this term}, and so the methods do not work
directly. However, we can get around this with a trick.

Note that any reciprocal polynomial in $x$ of even degree
$2n$ can be written (uniquely) as a multiple (by $x^n$) of a
polynomial $g(y)$ in $y=x+1/x$ of degree $n.$ The proof is very
simple: Dividing through by $x^n,$ we write 
\[f(x) = a_{n} + \sum_{i=0}^{n-1} a_i(x^{i-n} + x^{n-i}).\]
Note that $(x+1/x)^n$ is a reciprocal polynomial, and so is
$f(x) - (x+1/x)^n,$ which is also of lower degree than $f(x).$ The
result now follows by induction (notice that the coefficients of $g$
are integer linear combinations (whose coefficients depend only on the
degree of $f$) of the coefficients of $f,$ and,
obviously, \emph{vice versa}.

Now, it is clear that in order for $f(x)$ to be reducible, $g(y)$ must
be also. Indeed, suppose $f(x) = f_1(x) \dots f_k(x),$ where the $f_i$
are irreducible. Since $x^2n f(1/x) = f(x),$ it follows that 
$f(x) = \prod_{i=1}^k x^{\deg f_i(x)} f_i(1/x).$ By the irreducibility
of $f_i(x),$ it follows that \emph{either} $f_i$ is a reciprocal
polynomial, or $f_i$ is the reciprocal of some $f_j,$ in which case
$f_i(x) f_j(x)$ is a reciprocal polynomial. So, $f(x)$ has a
reciprocal factorization, and so $g(y)$ is reducible.

We now reason as in Section \ref{arbitrarypoly}, but with polynomials
$g(y)$ replacing $f(x).$

\section{Applications to matrices}
\label{matrices}

\subsection{The special linear group.}
\label{slnsec}
Consider first the matrix group $\SL(n, k).$ Since the coefficients of the
characteristic polynomial of a matrix $M$ are polynomials in the
entries of $M,$ and the dimension of $\SL(n)$ is $n^2 - 1,$ we see that
\begin{lemma}
\label{slnpgro}
The number of matrices in $\SL(n, p)$ whose characteristic
polynomial has a factor over $F_p,$ with constant term $1,$
grows as $o(p^{n^2-1}).$
\end{lemma}
\begin{proof}
The proof requires one additional observation: that every monic
polynomial $p(x)$ of degree $d$ with constant term $1$ is the characteristic polynomial of
some matrix in $\SL(d)$ -- namely the \emph{companion matrix} of
$p(x).$  It follows that the set of matrices whose characteristic
polynomial satisfies the assumptions of the Lemma lies on an algebraic
subvariety of $\SL(n),$ and the result follows by Lang-Weil.
\end{proof}
\begin{corollary}
\label{slnpprob}
The probability that a matrix in $\SL(n, p)$ satisfies the hypotheses
of Lemma \ref{slnpgro} goes to $0$ as $p$ goes to infinity.
\end{corollary}
\begin{proof}
The order of $\SL(n, p)$ is well known to be 
\[
p^{(n^2-n)/2} (p^2 - 1) (p^3 - 1) \dots (p^n - 1) \sim p^{n^2-1},
\] (see Newman's book \cite{newmanmats}[VII.17]). The assertion of the
corollary follows immediately.
\end{proof}

Unfortunately, since the number of integral points on $\SL(n,
\mathbb{Z})$ of height (absolute value) bounded by $B$ grows much slower than
$B^{n^2-1}$ the above results do not imply the following
\begin{conjecture}
\label{slnzgro}
The probability that a matrix in $\SL(n, \mathbb{Z})$ with coefficients
bounded by $B$ has reducible characteristic polynomial goes to $0$ as
$B$ goes to infinity.\footnote{It has been suggested by Peter Sarnak
  that the methods of \cite{DRS} can be extended to prove this
  conjecture.This is the subject of a forthcoming paper by the author.}
\end{conjecture}
But since we know that the number of points on $M^{m\times n}(\mathbb{Z})$ of
height bounded by $B$ grows like $B^{n^2},$ we do have 
\begin{theorem}
\label{glnzgro}
The probability that a matrix in $M^{n\times n}(\mathbb{Z})$ with coefficients
bounded by $B$ has reducible characteristic polynomial goes to $0$ as
$B$ goes to infinity.
\end{theorem}
\begin{proof} The probability that such a matrix factors modulo a
large prime $B < p < 2B$ (factors having constant terms equal to the
divisors of the constant term of the characteristic polynomial mod
$p$) already goes to $0,$ 
\end{proof}

\subsection{Lower bounds and asymptotics}
Theorem \ref{glnzgro} gives an estimate of $O(B^{n^2-1}\log B)$ on the
number of matrices in $M^{n\times n}(\mathbb{Z})$ with reducible
characteristic polynomial. To get a \emph{lower} bound, we recall the
following theorem of Yonatan Katznelson:
\begin{theorem}[Y. Katznelson, \cite{katznelsonsing}]
\label{katzthm}
The number of $n \times n$ singular integral matrices with entries bounded by $B$ is
asymptotic to $c_n B^{n^2-n} \log B.$
\end{theorem}

The following Corollary is quite easy:
\begin{corollary}
\label{katzcor}
The number of $n\times n$ matrices whose characteristic polynomial has
a \emph{linear} factor over $\mathbb{Z}$ is bounded below by $c^\prime_n
B^{n^2-n+1}\log B.$
\end{corollary}
\begin{proof}
For every singular matrix $M,$ the matrices $M + k I_n,\quad k \in
\mathbb{Z}$ have characteristic polynomial which has a linear factor
over $\mathbb{Z}.$
\end{proof}

So, it follows that if $N_{n, B}$ is the number of reducible integer
matrices with coefficients bounded by $B,$ we have, for some non-zero
constants $c_1, c_2:$
\begin{equation}
\label{glbds}
c_1 B^{n^2-n +1} \log B \leq N(n, B) \leq c_2 B^{n^2-1} \log B.
\end{equation}
Note that for $n=2,$ the upper and lower bounds grow at the same rate,
so we now the order of growth (which can be sharpened to an asymptotic
result without too much difficulty). Otherwise, there is a
considerable gap between the upper and the lower bounds, We conjecture
that the lower bound is the truth:
\begin{conjecture}
\[
N(n, B) \asymp c_n B^{n^2-n+1} \log B.
\]
\end{conjecture}

\section{Random products of matrices in the symplectic and special
linear groups}
\label{matrices2}

In the preceeding section we defined the size of a matrix by (in
essence) its $L^1$ norm (any other Banach norm will give the same
results). However, it is sometimes more natural to measure size
differently: In particular, if we have a generating set $\gamma_1,
\dots, \gamma_l$ of our lattice $\Gamma$ (which might be $\SL(n,
\mathbb{Z})$ or $Sp(2 n, \mathbb{Z})$) we might want to measure the
size of an element by the length of the (shortest) word in $\gamma_i$
equal to that element -- this is the combinatorial measure of
size. The relationship between the size of elements and combinatorial
length is not at all clear, so the results in this section are proved
quite differently from the results in the preceding section. We will
need the following results: First a result of this author 
\begin{theorem}[Rivin \cite{walks}]
\label{mythm}
Let $G$ be a graph whose vertices are labeled by generators of a
\emph{finite} group $\Gamma.$ Consider the set of $S_N$ elements of $\Gamma$
obtained by multiplying elements along walks of length $n.$ Then,
$S_N$ becomes equidistributed over $\Gamma$ as $N$ goes to infinity.
\end{theorem}

We will also need the following results of Nick Chavdarov and Armand Borel.

\begin{theorem}[Chavdarov, A.~Borel \cite{chavdarov}]
\label{chavthm}
Let $q > 4$, and let $R_q(n)$ be the set of $2n \times 2n$ symplectic matrices over the field $F_q$ with
\emph{reducible} characteristic polynomials. Then 
\[
\dfrac{|R_q(n)|}{|\Sp(2n, F_p)|} < 1- \frac{1}{3 n}.
\]
\end{theorem}

\begin{theorem}[Chavdarov, A.~Borel \cite{chavdarov}]
\label{chavthm2}
Let $q > 4$, and let $G_q(n)$ be the set of $n \times n$ matrices with
determinant $\gamma \neq 0$ over the field $F_q$ with
\emph{reducible} characteristic polynomials. Then 
\[
\dfrac{|G_q(n)|}{|\SL(n, F_q)|} < 1- \frac{1}{2 n}.
\]
\end{theorem}

Theorem \ref{chavthm2} follows easily from the following result of A.~Borel:
\begin{theorem}[A.~Borel]
\label{borthm2}
Let $F$ be a monic polynomial of degree $N$ over $\mathbb{Z}/p \mathbb{Z}$
with nonzero constant term.  Then, the number $\#{F, p}$ of matrices in $\GL(N,
p)$ with characteristic polynomial equal to $F$ satisfies
\[
(p-3)^{N^2-N} \leq \#(F, p) \leq (p+3)^{N^2-N}.
\]
\end{theorem}
Theorem \ref{borthm2} will be used in Section \ref{strong}.
A result we will need in Section \ref{mcg}, and might as well state here,  is:

\begin{theorem}[D. Kirby, \cite{Kirby_symplecticchar}]
\label{kirby}
Any reciprocal polynomial is the characteristic polynomial of a symplectic
matrix. 
\end{theorem}

We now have our results:

\begin{theorem} 
\label{randprodthm}
Let $G$ and $S_N,$ be as in the statement of Theorem
\ref{mythm}, but with $\Gamma = \Sp(2n, \mathbb{Z}),$ or $\Gamma = \SL(2,
\mathbb{Z}).$ 
 Then the
probability that a matrix in $S_N$ has a reducible characteristic
polynomial goes to $0$ as $N$ tends to infinity.
\end{theorem}
\begin{proof}  Let $\Gamma_l$ be the set of matrices in $\Gamma$ reduced
  modulo $l$ -- it is known (see \cite{newmanmats}) that $Gamma_l$ is $\SL(n,
  l)$ or $\Sp(2n, l)$ (depending on which $\Gamma$ we took. Let $p_1, \dotsc,
  p_k$ be distinct primes, let $K = p_1 
\dots p_k.$ 
We know that:
\[\Gamma_K = \Gamma_{p_1} \times \dots
\times \Gamma_{p_k}.\] (see \cite{newmanmats} for the proof of the last
equality). A generating set of $\Sp(2n, \mathbb{Z})$ projects via
reduction modulo $K$ to a a generating set of $\Gamma_K$ (see, again,
Newman's book \cite{newmanmats}), and also, via reduction mod $p_i$ to
generating sets of the $\Sp(2n, p_i).$ By Theorems \ref{mythm} and
\ref{chavthm}, the probability that the characteristic polynomial in a
random product of $N \gg 1$ generators is reducible modulo all of the
$p_i$ is at most equal to $(1-3/n)^k.$ Since this is an upper bound on
the probability of being reducible modulo $\mathbb{Z},$ the result
follows.
\end{proof}

\begin{remark}
Using Lemma \ref{slnpgro} instead of Theorem \ref{chavthm2} for
$\SL()$ gives a sharper result, as well as a more elementary argument.
\end{remark}

An example of a graph $G$ is a 
 bouquet of circles. In this case, we are just taking
random products of generators or their inverses. Another is the graph (studied
 in \cite{walks}) where a generator is never followed by its inverse (so only
 reduced words in generators are allowed), and so on.

\section{Stronger irreducibility}
\label{strong}
We might ask if something stronger than irreducibility of the
characteristic polynomial can be shown.The answer is in the
affirmative. Indeed, the methods of the preceeding sections combined with the
results of the Appendix give immediately:
\begin{theorem}
\label{galthm}
The probability that a random word of length $L$ in a generating set of $SL(N,
\mathbb{Z})$ has characteristic polynomial with Galois group $S_N$
goes to $1$ as $L$ goes to infinity.
\end{theorem}

Aside from its intrinsic interest, Theorem \ref{galthm} implies the
following:
\begin{theorem}
\label{iwip}
The probability that a random word $w$ of length $L$ in a generating set
of $\SL(N, \mathbb{Z})$ and all proper powers $w^k$ have irreducible
characteristic polynomials goes to $1$ as $L$ goes to infinity.
\end{theorem}

Theorem \ref{iwip} will follow easily from Theorem \ref{galthm}
together with the following Lemma:
\begin{lemma}
\label{iwip2}
Let $M \in \SL(n, \mathbb{Z})$ be such that the characteristic
polynomial of $M^k$ is \emph{reducible} for some $k.$ Then the Galois group
of the characteristic polynomial of $M$ is \emph{imprimitive}, or the characteristic polynomial of $M$
is cyclotomic.
\end{lemma}
\begin{remark} For the definition of \emph{imprimitive} see, for
  example, \cite{wielandt,mhall}.
\end{remark}
\begin{proof}
Assume that the characteristic polynomial $\chi(M)$ is irreducible
(otherwise the conclusion of the Lemma obviously holds, since the
Galois group of $\chi(M)$ is not even transitive). Let the roots of
$\chi(M)$ (in the algebraic closure of $\mathbb{Q}$) be $\alpha_1, \dotsc, \alpha_n.$ The roots of
$\chi(M^k)$ are $\beta_1, \dotsc, \beta_n,$ where $\beta_j = \alpha_j^k.$ Suppose that
$\chi(M^k)$ is reducible, and so there is a factor of $\chi(M^k)$
whose roots are $\beta_1, \dots, \beta_l,$ for some $l < n.$ Since
$\Gal(\chi(M))$ acts transitively on $\alpha_1, \dotsc, \alpha_n,$ it
must be true that for every $i \in \{1, \dots, n\},$ $\alpha_i^k =
\beta_j,$ for some $j \in \{1, \dotsc, l\}.$ Let $B_j$ be those $i$
for which $\alpha_i^k = \beta_j.$ This defines a partition of $\{1,
\dotsc, n\}$ into blocks, which is stabilized by the Galois group of
$\chi(M),$ and so $G$ is an intransitive subgroup of $S_n,$
\emph{unless} $l = 1.$ In that case, the characteristic polynomial of
$M^k$ equals $(x-\beta)^n,$ and since $M^k \in \SL(n, \mathbb{Z})$ it
follows that $\beta = 1,$ and all the eigenvalues of $M$ are $n$-th
roots of unity, so that $M^k = 1.$
\end{proof}

\section{The mapping class group} 
\label{mcg}
Let $S_g$ be a closed surface of genus $g,$
 and let $\Gamma_g$ be the mapping class group of $S_g.$ The group $\Gamma_g$
 admits a homomorphism $\mathfrak{s}$ onto $\Sp(2g, \mathbb{Z})$  (we associate to each
 element its action on homology; the symplectic structure comes from the
 intersection pairing). The following result can be find in \cite{CassonBleiler}:
\begin{theorem}
\label{cassonthm}
For $\gamma \in \Gamma_g$ to be pseudo-Anosov, it is sufficient that $g =
\mathfrak{\gamma}$ satisfy all of the following conditions:
\begin{enumerate}
\item
The characteristic  polynomial of $g$ is irreducible.
\item
The characteristic polynomial of $g$ is not cyclotomic.
\item
The characteristic polynomial of $g$ is not of the form $g = h(x^k),$ for some
$k>1.$ 
\end{enumerate}
\end{theorem}
The following is a corollary of our results on matrix group:
\begin{theorem}
Let $g_1, \dots, g_k$ be a generating set of $\Sp(2n, \mathbb{Z}).$ The
probability that a random product of length $N$ of  $g_1, \dots, g_k$
satisfies the conditions of Theorem \ref{cassonthm} goes to $1$ as $N$ goes to
infinity. 
\end{theorem}
\begin{proof}
We prove that the probability that the random word $w_N$ \emph{not} satisfy the
conditions goes to $0.$ By Theorem \ref{randprodthm}, the probability that
$w_N$ has reducible characteristic polynomial goes to $0.$ In order for the
characteristic polynomial to be of the form $g = h(x^k)$ it is necessary that
the linear term (the trace) vanish. This is a proper subvariety of $\Sp(2g),$
and so the number of elements of any $\Sp(2g, p)$ satisfying this condition
is of order of $p^{2g^2 + g - 1}.$ Since the number of elements in $\Sp(2g,
p)$ is of order of $p^{2g^2 + g}$ (Dickson's Theorem, see \cite{newmanmats}),
the proof of Theorem \ref{randprodthm} goes through verbatim (but needs Theorem
\ref{kirby}) to show that this
is an asymptotically negligible condition. Finally, since the set of
cyclotomic polynomials of a given degree is finite, the set of symplectic
matrices with those characteristic polynomials is also a subvariety of the
full group (again, needing Theorem \ref{kirby}), and the same result holds.
\end{proof}

\section{Free Group Automorphisms}
\label{free}

An automorphism of $\phi$ of  a free group $F_n$ is called
\emph{strongly irreducible}\footnote{This terminology, with strong
  support from this author, has been introduced by L. Mosher and
  M. Handel for what was previously known as \emph{irreducible with irreducible powers}} if no (positive) power of
$\phi$ sends a free factor $H$ of $F_n$ to a conjugate. This concept was
introduced by M.~Bestvina and M.~Handel \cite{bestvinahandel1}, and many of the results of the
theory of automorphisms of free groups are shown for such
automorphisms. By passing to the action of $\phi$ on homology, Section
\ref{strong}\footnote{We need to change $\SL(n, \mathbb{Z})$ to
  $\GL(n, \mathbb{Z})$ throughout} shows the following:

\begin{theorem}
\label{irredpow}
Let $f_1, \dots, f_k$ be a generating set of the automorphism group of
$F_n.$ Consider all words of length $L$ in $f_1, \dots, f_k.$ Then,
for any $n,$ the probability that such a word is irreducible tends to
$1$ as $L$ tends to infinity and also the probability that
such a word is strongly irreducible tends to $1$ as $L$
tends to infinity.
\end{theorem}
\appendix
\section{Galois groups of generic restricted polynomials}

Let $P_{N, d}(\mathbb{Z})$ be the set of monic polynomials of degree $d$ with
integral coefficients bounded by $N$ in absolute value. It is a
classical result of B.~L.~van~der~Waerden that the probability that
the Galois group of $p \in P_{N, d}(\mathbb{Z})$ is the full symmetric
group $S_d$ tends to $1$ as $N$ tends to infinity. The argument is
quite elegant: First, it is observed that a subgroup $H < S_d$ is the
full symmetric group if and only if $H$ intersects every conjugacy
class of $S_d.$ This means that $H$ has an element with every possible
cycle type. It is further noted that there is a cycle type $(n_1,
\dots, n_k)$ in the Galois group of $p$ over $\mathbb{Z}/p \mathbb{Z}$
if and only if $p$ factors over $\mathbb{Z}/p\mathbb{Z}$ into
irreducible polynomials of degrees $n_1, \dotsc, n_k.$ Using
Dedekind's generating function for the number of irreducible
polynomials over $\mathbb{Z}/p\mathbb{Z}$ of a given degree, it is
shown that the probability of a fixed partition is is bounded below by
a constant (independent of the prime $p$), and the proof is finished
by an application of a Chinese Remainder Theorem.

In this note, we ask the following simple-sounding question: Let $P_{N, d, a,
k}(\mathbb{Z})$   be the set of all polynomials in $P_{N,
d}(\mathbb{Z})$ where the coefficient of $x^k$ equals $a.$ Is it still
true that the Galois group of a random such polynomial is the full
symmetric group? The result would obviously follow if the probability
that the Galois group of a random general polynomial is ``generic''
were to go to $1$ sufficiently fast with $N.$ In fact, the probability
that an element of $P_{N, d}$ is \emph{reducible} (which means that
its Galois group is not transitive, hence not $S_n$) is of the order
of $1/N,$ so that approach does not work. 

Mimicking the proof of van der Waerden's result does not appear to work
(at least not easily): Dedekind's argument enumerates all irreducible
polynomials, and the result is not ``graded'' by specific
coefficients. It is certainly possible that the argument can be pushed
through, but this appears to be somewhat involved. 

Given this sad state of affairs, we first use a simple trick and Dirichlet's
theorem on primes on arithmetic progressions to show first the following technical result:
\begin{theorem}
\label{premainthm}
The probability that a random element of $P_{N, d, a, k}(\mathbb{Z}/p
\mathbb{Z})$ has a a prescribed splitting type $s$ 
approachs  the probability that a random unrestricted polynomial of degree
$d$ has the splitting type $s,$ as long as $p-1$ is relatively prime
to $(d-k)!,$ and as $p$ becomes large.
relatively prime to $d-k.$
\end{theorem}

which implies (by van der Waerden's sieve argument):

\begin{theorem}
\label{mainthm}
The probability that a random element of $P_{N, d, a, k}(\mathbb{Z})$
has $S_d$ as the Galois group tends to $1$ as $N$ tends to infinity,
\end{theorem}

It should be noted that the (multivariate) Large Sieve (as used by  P. X. Gallagher
in \cite{galgal}) can be used to give an effective estimate on the
probability in the statement of Theorems \ref{mainthm}: that
is: $p(N) \ll N^{-1/2}.$

\subsection{Proof of Theorem \ref{mainthm}}
\label{mainproof}
We will need two ingredients other than van der Waerden's original 
idea. The first of these is A.~Weil's estimate on the number of
$\mathbb{F}_p$ points on a curve defined over $\mathbb{F}_p:$
\begin{theorem}[A. Weil,\cite{weilcourbes}]
\label{weilcourbes}
Let $f \in \mathbb{F}_p[X, y]$ be an absolutely irreducible (that is,
irreducible in $\overline{\mathbb{F}}_p[X, Y]$) polynomial of degree $d$. Then if
\[\mathcal{C} = \{(x, y) \in \mathbb{F}_p^2 | f(x, y) = 0\},\]
we have the estimate
\[||C| - p| \leq 2 g \sqrt{p} + d^2,\] where $g$ is the genus of the
curve defined $f$ (which satisfies $g \leq (d-1)(d-2).$
\end{theorem}
This estimate is optimal.

The other classical result we shall need is the following:
\begin{theorem}[\cite{lang}[Theorem VIII.9.1]]
\label{langthm}
Let $k$ be a field, and $n\geq 2$ an integer. Let $a\in k,$ $a\neq 0.$
Assume that for all prime numbers $p$ such that $p|n$ we have $a
\notin k^p,$ and if $4|n,$ then $a\notin - 4 k^4.$ Then $X^n - a$ is
irreducible in $k[X].$
\end{theorem}

Theorem \ref{langthm} goes essentially back to N. H. Abel's
foundational memoir.

We will need an additional observation:

\begin{lemma}
\label{freelem}
Let $q = p^l,$ and 
let $x_1, \dots, x_k \in \mathbb{F}_q.$ Let $a, b \in \mathbb{F}_p,$
with $(a, b) \neq (1, 0).$
and let $g(a, b)(x) = a x + b$ be a transformation of $\mathbb{F}_q$
to itself. Then, it is not possible for $g(a, b)$ to permute $x_1,
\dots, x_k,$ if $k!$ is coprime to $p-1.$
\end{lemma}

\begin{lemma}
\label{freecor}
Consider a polynomial $f$ of degree $d$ over $\mathbb{F}_p,$ such that
$d < p,$ and such that the coefficient of $x^{d-1}$ does not vanish. Then there is no pair $(a, b) \neq (1, 0),$ such that
$f(a x + b) = a^d f(x),$ for all $x \in \mathbb{F}_p.$
\end{lemma}
\begin{proof}
There are two distinct cases to analyze. The first is when $a=1.$ In
that case, $f(x+b) = f(x)$ for all $x\in \mathbb{F}_p,$ and since $p >
d,$ $f(x+b) = f(x),$ for all $x$ in the algebraic closure of $F_p.$ Let
$r$ be a root of $f.$ Then, so are $r+a, r+2a, \dotsc, r+a(p-1),$ but
since $p$ is greater than $d$ that means that $f$ is identically $0.$

The second case is when $a\neq 1.$ In that case, $x_0 = b/(1-a)$ is
fixed under the substitution $x \rightarrow a x + b,$ and changing
of variables to $z = x - x_0,$ sends $f(z)$ to $f(a z).$ By the same
argument as above, $f(az) = a^d f(z),$ and so the corresponding
coefficients of the right and the left hand polynomials must be equal  Since the coefficient of
$x^{d-1}$ does not vanish, it follows tha $a=1,$ which contradicts our
assumption. 
\end{proof}

The argument now proceeds as follows. First, we note that if the
polynomial $f(x)$ of degree $d$ has a certain splitting type (hence Galois group)
over $\mathbb{F}_p$ then so does $f(a x + b)/a^d,$ for any $a\neq 0, b \in
\mathbb{F}_p.$ The set of all linear substitutions forms a group $\mathbb{A}$,
which acts freely on the set of polynomials of degree $d,$ except for
the (small) exceptional set of polynomials with a vanishing
coefficient of $x^{d-1}$ as long as $d < p$ (by Lemma \ref{freecor}), so the distribution of splitting types among the $\mathbb{A}$
orbits is the same as among all of the polynomials of degree $d.$ Now,
consider polynomials with constant term $1.$ How many of them are
there in the $\mathbb{A}$ orbit of $f(x)?$ It is easy to see that the
number is equal to the number of solutions to 
\[f(b) = a^d.\] If the curve $\mathcal{C}_f$ given by $f(x) - y^d$ is absolutely irreducible,
that number is $p + O(\sqrt{p}),$ by Theorem \ref{weilcourbes}. By
Theorem \ref{langthm}, in order for $\mathcal{C}_f$ to not be
absolutely irreducible, we must either have that $f(x) = g^q(x),$ for
some $q|d,$ or $f(x) = - 4 h^4(x),$ in case $4|d.$ But the number of
such polynomials is bounded by $O(p^{d/2}),$ which is asymptotically
neglible. So, we see that the distribution of splitting types amongst
polynomials of degree $d$ with constant term $1$ is the same as for
all polynomials, as long as $d < p.$ 

\bibliographystyle{plain}
\bibliography{rivin}
\end{document}